\input amssym.def
\input amssym
\magnification=\magstep1
\def\nbigskip{\bigskip\noindent}
\def\bbigskip{\bigskip \bigskip  }
\def\nbbigskip{\bbigskip \noindent }
\def\nbbbigskip{\bigskip\bigskip\bigskip\noindent}
\def\nmedskip{\medskip\noindent}
\def\qq{/\kern-.185em /}
\def\C{{\Bbb C}}
\def\Z{{\Bbb Z}}
\def\N{{\Bbb N}}
\def\R{{\Bbb R}}
\def\P{{\Bbb P}}
\def\Q{{\Bbb Q}}
\def\im{{\rm Im\,}}

\def\Gl{{\rm GL}}

\def\H{{\rm H}}
\def\buildunder#1#2{\mathrel{\mathop{\kern0pt #2}\limits_{#1}}}
\def\buildover#1#2{\buildrel#1\over#2}
\def\qed{\hfill{$\square$}}
\def\Lie{{\rm Lie}}

\def\psh{plurisubharmonic}
\def\spsh{strictly \psh}
\def\litem{\par\noindent\hangindent=\parindent\ltextindent}
\def\ltextindent#1{\hbox to \hangindent{#1\hss}\ignorespaces}
\def\litem{\par\noindent\hangindent=\parindent\ltextindent}
\def\ltextindent#1{\hbox to \hangindent{#1\hss}\ignorespaces}
\def\litem{\par\noindent\hangindent=\parindent\ltextindent}
\def\ltextindent#1{\hbox to \hangindent{#1\hss}\ignorespaces}
\def\litem{\par\noindent\hangindent=\parindent\ltextindent}
\def\ltextindent#1{\hbox to \hangindent{#1\hss}\ignorespaces}

\centerline{\bf  Projectivity of moment map quotients}

\bigskip
\centerline{\sl Peter Heinzner and Luca Migliorini}

\bigskip\bigskip\bigskip
Let $G$ be a complex reductive group acting algebraically on 
a complex projective variety $X$. Given a polarisation of $X$,
i.e., an ample $G$-line bundle $L$ over $X$, Mumford (see [M-F-K])
defined the notion of stability: A point $x\in X$ is
said to be semistable with respect to $L$ if and only if there
exits $m\in \N$ and an invariant section 
$s:X\to L^m$ such that $s(x)\ne 0$. Let $X(L)$ denote the set of 
semistable points in $X$, then there is a projective variety 
$X(L)\qq G$ and a $G$-invariant
surjective algebraic map $\pi:X(L)\to X(L)\qq G$ such that 

\smallskip
\item{(i)} $\pi$ is an affine map and

\smallskip
\item{(ii)} ${\cal O}_{X(L)\qq G}=(\pi_*{\cal O}_{X(L)})^G$.

\smallskip
In particular, for an open affine subset $U$ of $X(L)\qq G$, it follows 
that $\pi^{-1}(U)=\hbox{Spec}\,\C[U]^G$ where $\C[U]$ denotes the coordinate
ring of $\pi^{-1}(U)$ and $\C[U]^G$ is the algebra of invariant
functions.

\smallskip
There is a completely analogous picture for a holomorphic action
of a complex reductive group $G$ on a K\"ahlerian space $X$. The 
role of a polarisation is taken over by a Hamiltonian action 
of a maximal compact subgroup $K$ of $G$. Here one considers
a maximal compact subgroup $K$ of $G$, assumes the
K\"ahler structure to be  $K$-invariant and that there is 
an equivariant moment map $\mu:X\to \frak k^*$ with respect to $\omega$. 
In this situation 
$X(\mu)=\{x\in X;\ \overline{G\cdot x}\cap \mu^{-1}(0)\ne\emptyset\}$
is called the set of semistable points of $X$ with respect to $\mu$.
Here $\overline{G\cdot x}$ denotes the topological closure of
the $G$-orbit through $x$. The following result has been proved in 
[H-L] (c.f. [S]).

\smallskip
The set $X(\mu)$ is open in $X$ and there is a complex space $X(\mu)\qq G$ 
and a $G$-invariant surjective holomorphic map $\pi:X(\mu)\to X(\mu)\qq G$
such that 

\smallskip
\item{(i)} $\pi$ is a Stein map and

\smallskip
\item{(ii)} ${\cal O}_{X(\mu)\qq G}=(\pi_*{\cal O}_{X(\mu)})^G$.

\smallskip
In fact there is one more analogy between these two constructions. In the 
case where $X$ is projective, the line bundle $L$ induces a line bundle 
$\bar L$ on $X(L)\qq G$ which turns out to be ample. In the K\"ahler case 
$\omega$ induces a K\"ahlerian structure $\bar \omega$ on $X(\mu)\qq G$.

\smallskip
 A very ample $G$-line bundle $L$ over $X$ induces a 
$G$-equivariant holomorphic embedding of $X$ into $\P(V)$ where $V$ is
the dual vector space of the space of sections $\Gamma(X,L)$ and the
$G$-action on $\P(V)$ is induced by the natural linear $G$-action 
on $\Gamma(X,L)$. Now one may assume the $K$-representation to be unitary
and therefore the pull back of the Fubini-Study form $\omega_{\P(V)}$
to $X$ is a $K$-invariant K\"ahler form $\omega$ and the pull back 
of the natural moment map $\mu_{\P(V)}$ to $X$ gives a moment map
$\mu:X\to \frak k^*$. In this case, using a result of Kempf-Ness 
(see [K-N]), one checks that $X(\mu)=X(L)$, i.e., the set of 
Mumford-semistable subsets of $X$ is a subset of the set of 
momentum-semistable sets (see [K], [N] or sec. 3.).

\smallskip
Of course in general a given invariant K\"ahler form $\omega$ on a 
projective $G$-manifold $X$ may not be integral.
Therefore assoziated moment maps  are not in an obvous way related
to $G$-line bundles. Nevertheless, our goal here is to prove the 
following

\nbbigskip
{\bf  Semistability Theorem.} {\it Let $X$ be a smooth projective variety
endowed with a holomorphic action of a complex reductive group
$G=K^\C$, $\omega$ a $K$-invariant K\"ahler 
form and $\mu:X\to \frak k^*$ a $K$-equivariant moment map. Then there is a 
very ample $G$-line bundle $L$ over $X$ such that 
$$X(\mu)=X(L)\,.$$ }

\nbbigskip
Recently there has been some interest in the question of how $X(L)$ and 
$X(L)\qq G$ vary in dependence of $L$ (see e.g. [D-W], [T]). The above 
obviously implies that these results extend to the case where $\mu$
is moving.

\nbbbigskip
{\bf 1. Mumford quotients}

\nbbigskip
Let $G$ be a complex reductive group und $V$ a 
$G$-representation, i.e., there is given a holomorphic
homomorphism $\rho:G\to \Gl(V)$. Since $G$ is reductive, it is
in fact a linear algebraic group and $\rho$ is an algebraic map
(see e.g. [Ch]).
Moreover the algebra $\C[V]^G$ of $G$-invariant polynomials
is finitely generated. The corresponding affine 
variety  is denoted by $V\qq G$. The inclusion 
$\C[V]^G\hookrightarrow\C[V]$  induces a polynomial map $\pi:V\to V\qq G$
which turns out to be surjective. Explicitely $\pi:V\to V\qq G$ can be 
realized as follows. Let $q_1,\ldots, q_k$ be a set of generators of
the algebra $\C[V]^V$ and $q:=(q_1,\ldots,q_k)$. Then $Y:=q(V)$ is a 
Zariski-closed subset of $\C^k$ which is isomorphic with $V\qq G$. Under this 
isomorphism $\pi:V\to V\qq G$ is given by $q$. 

\smallskip
Since the group $G$ and the action $G\times V\to V,\, (g,v)\to g\cdot v,$
are algebraic, every $G$-orbit is Zariski-open in its closure. In particular,
for every $x\in\overline{G\cdot v}\setminus G\cdot v$ we have $\dim G\cdot x
<\dim G\cdot v$. This implies that the closure of every $G$-orbit 
contains  a closed $G$-orbit which may be defined as a $G$-orbit of smallest
dimension in $\overline{G\cdot v}$. Now $G$-invariant polynomials separate
$G$-invariant Zariski-closed subsets. This can be seen by 
using integration over
a maximal compact subgroup $K$ of $G$.  Thus the closed $G$-orbit in 
$\overline{G\cdot v}$ is unique. Moreover for $v,w\in V$ we have 
$\pi(v)=\pi(w)$ if and only if 
$\overline{G\cdot v}\cap \overline{G\cdot w}\ne\emptyset$ and this is the case 
if and only if $\overline{G\cdot v}$ and $\overline{G\cdot w}$ contain
the same closed orbit. Consequently, if $G\cdot v_0$ is the closed 
orbit in $\overline {G\cdot v}$, then 
$\pi^{-1}(\pi(v))=\{w\in V;\, G\cdot v_0\subset\overline{G\cdot w}\}$.
This is often expressed by the phrase that the quotient $V\qq G$ parametrises
the closed $G$-orbits in $V$.

\smallskip
Assume now that $X$ is a projective $G$-variety which is realized as a 
$G$-stable Zariski-closed subset of $\P(V)$. In general there is no 
way to associate
to $X$ a quotient $X\qq G$ which has reasonable properties. For 
example if $V$ is irreducible, then $\P(V)$ contains a unique $G$-orbit 
which is compact. This orbit is the image of a $G$-orbit through a maximal
weight-vector in $V$. Since every $G$-orbit in $\P(V)$ contains a 
closed $G$-orbit in its closure, the unique compact orbit is 
contained in the closure of every other $G$-orbit in $\P(V)$. If one 
would try to define a Hausdorff quotient, then every point would have to
be identified with the points in the unique compact
orbit. The resulting quotient would be a point. 

\smallskip
In order to resolve this difficulty
Mumford introduced  the following procedure (see [M-F-K]). Let $N$ be the 
Null-cone in $V$, i.e., 
the fibre through the origin of the quotient map $\pi:V\to V\qq G$ and
let $p:V\setminus \{0\}\to\P(V)$ denote the $\C^*$-principal bundle
which defines the projective space $\P(V)$.  For a subset $Y$ of $\P(V)$ let
$\hat Y:=p^{-1}(Y)$ denote the corresponding cone in $V$. A point 
$x\in X$ is said to be semistable with respect to $V$ if 
$\hat x=p^{-1}(x)\subset \hat X\setminus N$. Let 
$X(V):=p(\hat X\setminus N)$ denote the set of semistable points in 
$X$ with respect to the representation $V$. Thus $X(V)$ is obtained by removing the 
image of the Null-cone from $X$.   

\smallskip
The cone $C(X):=\hat X\cup \{0\}$ in $V$ over $X$ is a $G$-stable closed affine
subset of $V$ and $N$ is saturated with respect to $\pi_V:V\to V\qq G$. Thus 
$\hat X(V):=\hat X\setminus N=C(X)\setminus N$ is saturated with respect to 
$\pi_{\hat X}:\hat X\to \hat X\qq G$. In particular there is  a quotient 
$\hat\pi:\hat X(V)\to\hat X(V)\qq G$ which is given by restricting 
$\pi_V:V\to V\qq G$ to $\hat X(V)$. The $\C^*$-action on $V$ defined by 
multiplication commutes with the $G$-action and stabilises $\hat X(V)$. 
Thus there is an induced $\C^*$-action
on $\hat X(V)\qq G$ which can be described explicitly as follows.
Let $q_1,\ldots , q_k$ be a set of homogeneous generators of $\C[V]^G$ with
$\deg q_j=d_j$. The map $q:V\to \C^k$ is equivariant with respect to $\C^*$.
More precisely we have $q(t\cdot v)=(t^{d_1}q_1(v),\ldots, t^{d_k}q_k(v))$. 
Moreover $q(V\setminus N)=q(V)\setminus\{0\}\subset \C^k\setminus \{0\}$. 
Note that $\C^*$ acts properly on $\C^k\setminus \{0\}$. In particular
there is a geometrical quotient $\C^k\setminus \{0\}/\C^*=:\P(d_1,\dots,d_k)$ 
which is a projective variety.
This implies that $X(V)\qq G:=(\hat X(V)\qq G)/\C^*$ is also a projective 
variety, since it is a Zariski-closed subspace of $\P(d_1,\dots,d_k)$. The map
$\hat X(V)\to X(V)\qq G$ is $\C^*$-invariant and induces therefore an algebraic
map $\pi:X(V)\to X(V)\qq G$ which is the quotient map for the $G$-action on
$X(V)$.

\smallskip
There is a standard procedure to realize a given $G$-variety $X$ as a 
$G$-stable subvariety of some projective space $\P(V)$ where 
$V$ is a $G$-representation. For this assume that $L$ is a very ample 
line bundle over $X$ and let $\Gamma (X,L)$ denote the space of holomorphic
sections of $L$. Thus the natural map $\imath_L:X\to \P(V)$ which is given 
by evaluation where
$V:=\Gamma(X,L)^*$ is the dual of $\Gamma(X,L)$ is an embedding. 
Now if the $G$-action on $X$ lifts
to a $G$-action on $L$, then $V$ is a $G$-representation in a natural way
and $\imath_L$ is $G$-equivariant. The set $X(L):=\{x\in X;\ s(x)\ne 0
\ \hbox{for some invariant section  } s\in \Gamma(X,L^m),\ m\in \N\}$ coincides with $X(V)$ after
identifying $X$ with $\imath_L(X)\subset \P(V)$ and is called the set of 
semistable points of $X$ with respect to the $G$-line bundle $L$. Note that 
$X(L)$ depends on $L$ and on the lifting of the $G$-action to $L$.

\smallskip
The following two elementary facts concerning $G$-actions on line bundles
are often useful.

\nbbigskip
{\bf Lemma.} {\it Let $X$ be a connected projective $G$-variety.  
\item{(i)} If $L$ is ample, then there is a lifting of the $G$-action
to some positive power $L^m$ of $L$. 
\smallskip
\item{(ii)} Two liftings of the $G$-action to $L$ differ by a character of 
$G$.}

\nbigskip
{\it Proof.}
The last statement follows since $X$ is compact and therefore a $G$-action 
on the trivial bundle $X\times \C=L\otimes L^{-1}$ is given by 
$g\cdot(x,z)=(g\cdot x, \chi(g)z)$ where $\chi:G\to \C^*$ is a character of
$G$.

\smallskip
For the first statement one may assume that $G$ acts effectively. Since $G$
is connected the induced action on Pic$(X)$ is trivial. This implies 
that there is a subgroup $\tilde G$ of the automorphism group of $L$ and 
an exact sequence of the form 
$$1\to \C^*\to \tilde G\to G\buildover \alpha\to 1$$
where $\alpha$ is given by restricting $\tilde g\in \tilde G$ to the zero 
section $X\hookrightarrow L$. This sequence splits after replacing 
$G$ by a finite covering. Hence the $G$-action on $X$ lifts to 
$L^m$ for some positive $m$.
\qed

\nbbbigskip
{\bf 2. Moment map quotients}

\nbbigskip
Let $G$ be a complex reductive group which acts holomorphically on a complex
manifold $X$. Now choose a maximal compact subgroup $K$ of
$G$ and let $\omega$ be a $K$-invariant K\"ahler form on $X$. 
By definition the $K$-action on $X$ is said to be Hamiltonian with moment
map $\mu$ if there is given an equivariant smooth map $\mu$ from $X$ 
into the dual $\frak k^*$ of the Lie algebra $\frak k$ of $K$ such that
$$d\mu_\xi=\imath_{\xi_X}\omega \leqno(*)$$
for all $\xi \in \frak k$. Here $\xi_X$ denotes the vector field on
$X$ associated with $\xi$, $\mu_\xi=\mu(\xi)$ and $\imath_{\xi_X}\omega$
is the one form $\eta\to \omega(\xi_X,\eta)$. Note that for a connected 
manifold $X$ an equivariant moment map is 
uniquely defined by $(*)$ up to a constant in $\frak k^*$ which lies in
the set of fixed points. In particular, if the group $K$ is semisimple
then an equivariant moment map is unique. Moreover in the semisimple case 
it can be shown that $\mu$ exists for a given $K$-invariant K\"ahler form 
$\omega$ (see e.g. [G-S]) 

\nbigskip
{\it Example.} Let $\rho:X\to \R$ be a 
smooth $K$-invariant function, $\omega:=2i\partial\bar\partial\,\rho$
and let $\mu:X\to \frak k^*$ be the associated $K$-equivariant map which is
defined by $\mu_\xi=d\rho(J\xi_X)$. Here $J$ denotes the complex
structure tensor on $X$. A direct calculation shows that 
$d\mu_\xi=\imath_{\xi_X}\omega$ holds for every $\xi\in \frak k$. 
In particular, if $\rho$ is strictly plurisubharmonic, i.e., $\omega$
is K\"ahler, then $\mu$ is a moment map. We refer to 
$\mu=:\mu^\rho$ as the moment map given by $\rho$.
 
\nbigskip
Similar to the case of an ample $G$-line bundle there is a notion of
semistability with respect to $\mu$. A point $x\in X$ is said to 
be semistable with respect to $\mu$ 
if $\overline{G\cdot x}\cap \mu^{-1}(0)\ne\emptyset$. Let $X(\mu)$ 
denote the set of semistable points with respect to $\mu$.

\smallskip
The following is proved in [H-L] (see also [S]).

\nbbigskip
{\bf Theorem 1.} {\it The set of semistable points $X(\mu)$ is open in 
$X$ and the semistable quotient $\pi:X(\mu)\to X(\mu)\qq G$ exists.
The inclusion $\mu^{-1}(0)\hookrightarrow X(\mu)$ induces a homeomorphism
$\mu^{-1}(0)/K\cong X(\mu)\qq G$.}
\qed 

\nbbigskip
By a semistable quotient of a complex space $Z$ (see [H-M-P] for
more details) endowed with
a holomorphic action of $G$ we mean 
a complex space $Z\qq G$ together with a $G$-invariant surjective
map $\pi:Z\to Z\qq G$ such that:

\smallskip
\item{{\rm (i)}} The structure sheaf ${\cal O}_{Z\qq G}$ is given by
$(\pi_*{\cal O}_Z)^G$, i.e., the holomorphic functions on an open subset 
of $Z\qq G$ are exactly the invariant holomorphic functions on its inverse
image in $Z$. 

\smallskip
\item{{\rm(ii)}} The map $\pi:Z\to Z\qq G$ is a Stein map, i.e., the inverse
image of a Stein subspace of $Z\qq G$ is a Stein subspace of $Z$.

\nbbigskip
In [H-H-L] it is shown that 
each point $q\in X(\mu)\qq G$ has an open neighborhood $Q$ such
that $\omega=2i\partial\bar\partial\rho$ on $\pi^{-1}(Q)$ for some
$K$-invariant smooth function $\rho$. Furthermore, the moment map 
$\mu$ restricted to $\pi^{-1}(Q)$ is given by $\rho$, i.e., $\mu=\mu^\rho$. 
A result of Azad and Loeb (see [A-L]) asserts that 
$\rho$ is an exhaustion on $G\cdot x$ if $x\in \mu^{-1}(0)$, 
i.e., is bounded from below and proper. In particular $G\cdot x$ is 
closed in $X(\mu)$ for every $x\in \mu^{-1}(0)$. The converse is also true
in the following sense. If $G\cdot x$ is closed in $X(\mu)$, then 
$\mu(g\cdot x)=0$ for some $g\in G$.
Furthermore in [H-H] it is shown that the restriction of $\rho$ to each 
fibre over $Q$ is an exhaustion, i.e., is bounded from below and proper.
This Exhaustion Lemma  and also a refinement of it (see sec. 6)
will be used several times in the 
remainder of this paper. For example, it implies the following 
(see [H-H]).

\nbbigskip
{\bf Theorem 2.} {\it Let $X$ be a compact complex manifold with a 
holomorphic $G$-action and let $\mu:X\to\frak k^*$ be a moment
map with respect to a $K$-invariant K\"ahler form $\omega$. Let 
$\tilde\omega$ be a $K$-invariant K\"ahler form on $X$ which lies in
the cohomology class of $\omega$. Then there exists a moment map
$\tilde\mu:X\to\frak k^*$ with respect to $\tilde\omega$ such that
$$X(\mu)=X(\tilde\mu)\ .$$}

\nbigskip
{\it Proof.} We recall the argument given in [H-H].
Since $\tilde\omega$ is cohomologous to $\omega$ and $X$ is a compact 
K\"ahler manifold, there exists
a differentiable $K$-invariant function $f:X\to\R$ so that 
$\tilde\omega=\omega+2i\partial\bar\partial f$. Define 
$\mu^f:X\to\frak k^*$ by $\mu^f_\xi=J\xi_X(f)$ for 
$\xi\in\Lie\,K$ and set $\tilde\mu=\mu+\mu^f$. Then $\tilde\mu$ is
a moment map with respect to $\tilde\omega$. For every $x\in X(\mu)$
there exists a \spsh\ $K$-invariant function $\rho:Z\to\R$, where
$Z:=\overline{G\cdot x}\cap X(\mu)$, so that $\mu\vert Z=\mu^\rho$,
where $\mu^\rho$ is the moment map associated to $\rho$
(see [H-H-L]). Since $Z\cap\mu^{-1}(0)\not=\emptyset$, the above mentioned
Exhaustion Lemma  implies that $\rho:Z\to\R$ is an exhaustion.
Now $f$ attains its minimum and maximum on $X$ and $\rho$ is an 
exhaustion. Hence the \spsh\ $K$-invariant function $\tilde\rho:=\rho+f$
is also an exhaustion on $Z$. This shows that $Z\subset X(\tilde\mu)$,
i.e., $X(\mu)\subset X(\tilde\mu)$. By symmetry we 
have $X(\mu)=X(\tilde\mu)$.
\qed

\nbbigskip
If $G$ is a connected semisimple Lie group, then a 
moment map with respect to a $K$-invariant K\"ahler $\omega$ 
always exists and is unique. Thus in this case Theorem 2 shows 
that $X(\mu)$ depends only on the cohomology class of $\omega$.

\nbbbigskip
{\bf 3. Moment maps associated to representations}

\nbbigskip
Let $V$ be a $G$-representation whose restriction to the maximal compact 
subgroup
$K$ of $G$ is unitary with Hermitian inner product $<\ ,\ >$. Then $\rho:V\to \R$,
$\rho(z)={1\over 2}\vert\vert z\vert\vert^2={1\over 2}<z,z>$, is a $K$-invariant 
\spsh\ exhaustion  function on $V$ and consequently $V=V(\mu)$ where 
$\mu:V\to \frak k^*$, $\mu_\xi(z)={1\over i}<\xi z, z>=d\rho(J\xi z)$ is the 
moment map associated with $\rho$. The K\"ahler form 
$\omega_V=2i\partial\bar\partial\,\rho$ is given by $\omega_V(v,w)=-\im <v,w>$. 
Since in this case the restriction of $\rho$ to every $\pi$-fibre is an exhaustion,
we have $V(\mu)=V$ and the inclusion $\mu^{-1}(0)\hookrightarrow V$ induces a 
homeomorphism $\mu^{-1}(0)\cong V\qq G$ (see sec. 2 Theorem 1). The essential 
part of this statement has already been proved in [K-N]. 

\smallskip
Let $S:=S(V):=\{z\in V;\ \vert \vert z\vert \vert=1\}$ denote the unit sphere
in $V$. Note that $S$ is a co-isotropic submanifold of $V$ with respect to $\omega_V$,
i.e., $(T_zS)^{\perp_{\omega_V}}=T_z(S^1\cdot z)\subset T_zS$ where the circle group
$S^1=\{\lambda\in \C;\ \vert \lambda \vert =1\}$ acts on $V$ by multiplication. 
This is easily seen by using the orthogonal decomposition $T_zV=T_z(\C^* z)\oplus
W$ where $W:=T_zS\cap iT_zS$ denotes the complex tangent space of $S$ at $z$.
The complex structure on $W$ induces the standard complex structure on 
$\P(V)=S(V)/S^1$. Moreover since $S$ is co-isotropic, there is a unique symplectic 
structure $\omega_{\P(V)}$ on $\P(V)$ such that 
$\imath^*_Sp^*\omega_{\P(V)}=\imath_S^*\omega_V$. Here $p:(V\setminus \{0\})\to 
(V\setminus \{0\})/\C^*=\P(V)$ denotes the quotient map and 
$\imath_S:S\hookrightarrow V$ is the inclusion. Furthermore, the definition
of the complex structure and of $\omega_{\P(V)}$ are compatible so that 
$\omega_{\P(V)}$ is in fact a K\"ahler form on $\P(V)$. Up to a positive 
constant it is the unique K\"ahler form on $\P(V)$ which is invariant 
with respect to the unitary group ${\rm U}(V)$. Note that $\omega_{\P(V)}$
is determined by 
$p^*\omega_{\P(V)}=2i\partial \bar \partial \log \rho=
2i(-{1\over \rho^2}\partial\rho\wedge\bar \partial
\rho+{1\over \rho}\partial\bar\partial\rho)$.

\smallskip
The induced $K$-action on $\P(V)$ is again Hamiltonian. The moment map
is given by $(\mu_{\P(V)})_\xi([z])={2\over i}{<\xi z,z>\over \vert\vert z\vert\vert^2}
=d\log\, \rho(z)(J\xi z)$. In particular we have
$\overline{G\cdot [z]}\cap \mu_{\P(V)}^{-1}(0)\ne\emptyset$ if and only if
$\overline{G\cdot z}\cap \mu^{-1}(0)\ne\emptyset$ and this is the case 
if and only if $f(z)\ne 0$ for some $G$-invariant homogeneous polynomial 
$f$ on $V$.

\smallskip
Now let $X$ be a $G$-stable subvariety of $\P(V)$. The pull back of $\omega_{\P(V)}$
to $X$ induces a K\"ahlerian structure $\omega$ on $X$ und the $K$-action
is Hamiltonian with moment map $\mu:X\to \frak k^*$, $\mu=\mu_{\P(V)}\vert X$.
We call $\mu$ the standard moment map induced by the embedding into $\P(V)$. The 
above  construction shows the following well known

\nbbigskip
{\bf Lemma.} {\it Let $L$ be a very ample $G$-line bundle over $X$ and consider
$X$ as a $G$-stable subvariety on $\P(V)$ where the embedding is given by $\Gamma(X,L)$
and $V=\Gamma(X,L)^*$.
Then 
$$X(\mu)=X(L)\,,$$ 
i.e., the semistable points with respect 
to the standard moment map on $\P(V)$ are the semistable points 
with respect to $L$.}
\qed

\nbbbigskip
{\bf 4. The main result}

\nbbigskip
Let $G$ be a connected complex reductive group and $K$ a 
maximal compact subgroup of $G$, i.e., $G=K^\C$. By a $G$-variety
we mean in the following an algebraic variety together with an 
algebraic action of $G$.

\smallskip
Let $X$ be a smooth projective $G$-variety and $\omega$ a $K$-invariant 
K\"ahler form on $X$. Assume that the $K$-action is Hamiltonian 
with respect to $\omega$, i.e., there is a $K$-equivariant moment
map $\mu:X\to \frak k^*$, and denote by $X(\mu):=\{x\in X;\
\overline{G\cdot x}\cap \mu^{-1}(0)\ne\emptyset\}$ the set of semistable
points with respect to $\mu$.

\nbbigskip
{\bf Semistability Theorem.} {\it There is a very ample $G$-line bundle
$L$ over $X$ such that 
$$X(\mu)=X(L)\,.$$}

\nbbigskip
Here $X(L)$ denotes the set of semistable points in $X$ in the sense of 
Mumford, i.e., $X(L)=\{x\in X;\ s(x)\ne 0\ \hbox{for some}\
G\hbox{-invariant holomorphic section}\ s \ \hbox{of}\ L^m,\ m\in \N\}$.

\smallskip
The case where $\omega$ is assumed to be integral is well known and
follows rather directly from the definitions using standard Kempf-Ness
type arguments. In fact it is a consequence of Theorem 2 of sec. 2
and the Lemma in sec. 3. 

\smallskip
The proof in the general case is divided 
into two steps. In the first part we consider forms $\omega$ whose
cohomology class $[\omega]$ is contained in the $\R$-linear span
of the ample cone in $H^{1,1}(X)$. The second part of the proof is more 
involved. It is a reduction procedure to the first case.

\smallskip
At least implicitly (see e.g. [D-W], [K], [N]) the ample cone case 
seems to be known. 
In order to be complete we include a proof in the next paragraph.

\nbbbigskip
{\bf 5. The ample cone case}

\nbbigskip
In this section $G$ is a connected complex reductive group with a 
fixed maximal compact subgroup $K$ and $X$ is a smooth projective $G$-variety. 
Let $\omega$ be a $K$-invariant K\"ahler form and assume that
there is an equivariant moment map $\mu:X\to\frak k^*$. In this section 
we proof the following

\nbbigskip
{\bf Proposition.} {\it Assume that the cohomology class of $\omega$ lies
in the real linear span of the ample cone in $H^{1,1}(X)$. Then there 
exists a very ample $G$-line bundle $L$ over $X$ such that
$$X(\mu)=X(L)\,.$$}

\nbbigskip
{\it Proof.}
Since $X(\mu)$ essentially depends only on the cohomology class of
$\omega$ (see sec. 2 Theorem 2), we may assume that there are equivariant 
holomorphic embeddings $\imath_k: X\to \P(V_k)$, $k=1,\ldots, m,$ so that 
$$\omega=\sum a_k\,\imath^*_k\omega_{\P(V_k)}$$
where $a_k$ are positive real numbers.

\smallskip
Let $\imath: X\to \P(V_1)\times\cdots \times \P(V_m)$ be the diagonal 
embedding. Then 
$$\omega=\imath^* (\sum a_k\,\pi_k^*\omega_{\P(V_k)})\,,$$ 
where $\pi_k:\P(V_1)\times\cdots\times\P(V_m)\to \P(V_k)$ denotes the 
projection. Hence the moment map $\mu$ is the restriction of a moment map
on $\P(V_1)\times\cdots\times\P(V_m)$ with respect to 
$\sum a_k\, \pi^*_k\omega_{\P(V_k)}$ which also will be denoted by 
$\mu$. Since $X$ is closed in $\P(V_1)\times\cdots\times\P(V_m)$,
we have 
$$X(\mu)=(\P(V_1)\times\cdots\times\P(V_m))(\mu)\cap X\,.$$
Thus for the proof of the proposition we may assume that 
$X=\P(V_1)\times\cdots\times\P(V_m)$, 
$\omega=\sum a_k\,\pi_k^*\omega_{\P(V_k)}$ and the $G$-action is
given by a representation $G\to \Gl(V_1)\times\cdots\times\Gl(V_m)$.

\smallskip
Let $T$ be a maximal compact torus in $K$. Then $T^\C$ is a maximal
algebraic torus in $G$. We now reduce the proof of the proposition 
to the case where $G=T^\C$ as follows.  

\smallskip
Let $\mu_T:X\to \frak t^*$ be the moment map for the $T$-action which is 
induced by $\mu$ and the embedding $\frak t\hookrightarrow \frak k$. 
Then it follows that 
$$X(\mu)=\bigcap_{k\in K}k\cdot X(\mu_T)$$
by the Hilbert Lemma version in [K] (sec. 8.8.). Thus it is sufficient to 
show the following

\nbigskip
{\it Claim.} There exists a very ample $G$-line bundle $L$ over 
$X=\P(V_1)\times\cdots\times\P(V_m)$ such that 
$$X(\mu_T)=X(L_T)$$
where $X(L_T)$ denotes the set of semistable points with respect to $L$
if one considers $L$ as a $T^\C$-bundle.
 
\nbigskip 
The proposition follows from the above claim, since 
$$X(\mu)=\bigcap_{k\in K}k\cdot X(\mu_T)=\bigcap_{k\in K}k\cdot X(L_T)=X(L)\,.$$

In order to prove the claim one may proceed as follows.

\smallskip
Let $S=S_1\times\cdots\times S_m$ be the maximal torus in 
$\Gl(V_1)\times\cdots\times\Gl(V_m)$ which contains the image of $T$ 
and $\mu_k:\P(V_k)\to \frak s^*_k$ 
the standard moment map on $\P(V_k)$. We will consider $\mu_k$ as a
moment map with respect to $S=S_1\times\cdots\times S_m$ where the 
factors of $S$ different from $S_k$ act trivially on 
$X=\P(V_1)\times\cdots\times\P(V_m)$. Since 
$\omega=\sum a_k\,\omega_{\P(V_k)}$, the moment map 
$\mu:X\to \frak t^*$ is given by 
$$\mu=a_1\,\mu_1+\cdots+a_m\,\mu_m+c$$
where $c\in \frak t^*$ and $\mu_k$ now denotes the map from 
$X$ to $\frak t^*$ which is given by $\mu_k:X\to \frak s^*$ composed 
with the dual of $\frak t\to \frak s$. Now if $\tilde a_k$ and
are positive rational numbers and $\tilde c$ is rational, then 
$\tilde \mu:=\tilde a_1\,\mu_1+\cdots+\tilde a_m\,\mu_k+\tilde c$ is
a moment map with respect to 
$\tilde \omega:=\sum \tilde a_k\,\pi_k^*\omega_{\P(V_k)}$. Since 
$\tilde a_k$ and $\tilde c$ are rational, it follows that
there is a very ample $G$-line bundle $L$ over $X$ such that 
$X(L)=X(\tilde\mu)$. Thus we have to show the following

\bbigskip
{\it There exists $\tilde a_k$ and $\tilde c$ such that 
$X(\mu)=X(\tilde\mu)$.}
\bbigskip

This statement follows from  convexity properties of $\mu$ as follows.
Since $T$ is compact, the set $X^T$ of $T$-fixed points in $X$ is 
smooth. Let $X^T=\cup_{j\in J} F_j$ be the decomposition 
into connected components. Note that $\mu$ is constant on every 
$F_j$, $j\in J$. 
For a set $J$ let 
${\cal P}(J)$ be the set of subsets of $J$. 
We say that ${\cal L}\in {\cal P}(J)$ is $\mu$-semistable
if $0\in \hbox{Conv}\{\mu(F_j);\, j\in {\cal L}\}$ where Conv denotes 
the convex hull operation in $\frak t^*$. Let 
$X_{\cal L}:=\{x\in X;\, \overline{T^\C\cdot x}\cap F_j\ne\emptyset \ 
\hbox{for all } \ j\in {\cal L}\}$. Since $\mu(\overline {T^\C\cdot x})=
\hbox{Conv}\{\mu(F_j);\,\overline {T^\C\cdot x}\cap F_j\ne\emptyset\}$ 
(see [A]), it follows  that 
$$X(\mu)=\bigcup X_{\cal L}\,.$$ 
Here the union is taken over the elements
${\cal L}$ of ${\cal P}(J)$ which are $\mu$-semistable.
For a given $\mu$ denote by $I(\mu)$ the set of $\mu$-semistable 
subsets of $J$. We show now that if a collection of subsets is of the form 
$I(\mu)$, then $I(\mu)=I(\tilde\mu)$ for some positive rational 
$\tilde a_k$ and rational $\tilde c$. 

\smallskip
In order to see this, let $\Lambda_{kj}:=(\mu_k(F_j)\in \frak t^*$. Note that
$\Lambda_{kj}$ are integral points in $\frak t^*$. A subset
$I\subset {\cal P}(J)$ is of the form $I(\mu)$ if and only if there exist 
positive real numbers 
$a_k$ and 
$c\in \frak t^*$ such that for all ${\cal L}\in {\cal P}(J)$ the 
following holds.

\nbigskip
$0\in \hbox{ Conv}\ \{\sum_k a_k\Lambda_{kj} + c\,;\ j\in {\cal L}\}$ 
if and only if ${\cal L}\in I$.

\nbigskip
This condition is equivalent to a collection of linear inequalities with 
integral coefficients in the unknowns $a_k$'s and $c$ which have a real solution if and only if they have a rational one.
\qed

\nbbbigskip
{\bf 6. Cohomologous K\"ahler forms on orbits}

\nbbigskip
In this section let $G$ be a connected complex reductive group with 
maximal compact subgroup $K$ and let $X=G\cdot x_0$ be a $G$-homogeneous 
manifold. We assume that there are given $K$-invariant K\"ahler forms
$\omega^j$, $j=0,1,$ on $X$ which are cohomologous and set 
$$\omega^t=(1-t)\,\omega^0+t\,\omega^1\,, \ t\in [0,1]\,.$$
Moreover, assume that there are $K$-equivariant moment maps
$$\mu^t:X\to \frak k^*\,,\ t\in[0,1]$$
with respect to $\omega^t$ such that the dependence on $t$ is continuous.

\nbbigskip
{\it Remark.} We have 
$\mu^t=(1-t)\,\mu^0+t\,\mu^1+c^t$ where $c^t\in \frak z^*$
is a constant. Here $\frak z$ is the Lie algebra of the center of $K$. The 
goal of this section is to obtain some control about the semistable set
if $t$ varies.

\nbbigskip
{\bf Lemma.} {\it If $M_k^{t_0}\ne\emptyset$ for some $t_0\in [0,1]$, then
$\omega^t=2i\partial\bar\partial \rho^t$
where $\rho^t=(1-t)\,\rho^0+t\,\rho^1$ and $\rho^j:X\to \R$, $j=0,1$, are 
$K$-invariant smooth functions.}

\nbigskip
{\it Proof.}
Since $M^{t_0}_K\ne\emptyset$, the orbit $X=G\cdot x_0$ is 
a Stein manifold (see e.g. [H] or [H-H-L]). Now $\omega^0$ and $\omega^1$ 
are assumed to be cohomologous. Thus there is a $K$-invariant smooth function
$u:X\to \R$ such that $\omega^1-\omega^0=2i\partial\bar\partial\, u$.
On the other hand 
$\omega^{t_0}=2i\partial\bar\partial \,f$ for some $K$-invariant smooth 
function $f:X\to \R$ (see sec. 2 and [H-H-L]). Thus 
$\omega^t=2i\partial\bar\partial\, \rho^t$ where $\rho^0:=u-t_0f$
and $\rho^1:=u+(1-t_0)f$.
\qed

\nbbigskip
Now let $Z$ denote the connected component of the identity of the center of
$K$ and let $S$ be a semisimple factor of $K$. Thus $K=S\cdot Z$ and
$\frak k= \frak s\oplus \frak z$ on the level of Lie algebras. Let 
$\mu_S^t$ (resp. $\mu^t_Z$) be the moment map with respect to the $S$-action
(resp. $Z$-action), i.e., the composition of $\mu^t$ with
the dual of the inclusion $\frak s\hookrightarrow \frak k$ (resp.
$\frak z\hookrightarrow \frak k$). We also set 
$M^t_K:=(\mu^t)^{-1}(0)$, $M^t_S:=(\mu^t_S)^{-1}(0)$ and
$M_Z^t=(\mu^t_Z)^{-1}(0)$.

\nbbigskip
{\bf Proposition.} {\it If $M_k^{t_0}\ne\emptyset$ and if the set 
$X(\mu^t_Z)$ of $Z^\C$-semistable points is independent of $t$, then
there are pluriharmonic $K$-invariant functions $h^t:X\to \R$ which
depend continuously on $t$ such that 
$$\mu^t=\mu^{\rho^t+h^t}\,.$$ 
}

\nbigskip
{\it Proof.}
It follows from the definition of a moment map that it is unique up to 
a constant. Thus $\mu^{\rho^t}=\mu^t+c^t$ 
where $c^t$ is a $K$-invariant constant, i.e.,
$c^t\in \frak z^*$.

\smallskip
The proof of the Proposition will be reduced to the case
of a compact Abelian group $T\cong (S^1)^k$. In this situation we have 
$T^\C\cong (\C^*)^k$ and $\frak t=\Lie\,T\cong\R^k$. Moreover, for any 
$c^t\in\frak t^*$, $c^t=(c_1^t,\ldots ,c_k^t)$, the function 
$\tilde h^t(z_1,\ldots , z_k)=c_1^t \log\vert z_1\vert +\cdots +c^t_k
\log\vert z_k\vert$ 
is pluriharmonic on $T^\C$ and satisfies $\mu^{\tilde h^t}=c^t$.

\smallskip 
Let $x_0\in M^{t_0}_K$ and set $L:=K_{x_0}$. Then we
have the following orthogonal decomposition of the Lie algebra $\frak k$.
$$\frak k=\frak t\oplus\frak z_L\oplus\frak s$$
where   
$\frak z_L:=\frak z\cap(\frak s+\frak l),\ \frak z=\frak t\oplus\frak z_L$ 
and $\frak s+\frak l=\frak s\oplus \frak z_L$.

\smallskip
Note that $\frak z$ is the Lie algebra of the group $K/S$ and 
$\frak s+\frak l$ is the Lie algebra of the subgroup $S\cdot L$ of $K$.
Since $K$ is connected  
$K/SL=(K/S)/(SL/S)=:T$ is a compact connected Abelian group. Hence we have 
$T\cong(S^1)^k$ and 
$\Lie\,T\cong (\frak k/\frak s)/((\frak s+\frak l)/\frak s)\cong
\frak z/\frak z_L=\frak t$. 
Now identify $\frak k\cong\frak k^*$, i.e., we have the orthogonal splitting 
$$\frak k^*=\frak t^*\oplus\frak z_L^*\oplus\frak s^*\ .$$

\nmedskip
{\it Claim.} $c^t\in \frak t^*\ .$

\nmedskip
For the proof let $x_0\in M^{t_0}_K$ be given and note that $Z^\C\cdot x_0$ is
closed in $X(\mu^{t_0}_Z)=X(\mu^t_Z)$. Thus there are $x_t\in Z^\C\cdot x_0$
such that $\mu^t_Z(x_t)=0$. In particular we have
$c^t=\mu^{\rho^t}(x_t)$. 
Now let $\xi=\tau+ \lambda+\sigma$, where 
$\tau\in\frak t,\ \lambda\in\frak z_L$ and $\sigma\in\frak s$. 
Then, since the moment map is unique for a semisimple Lie group, 
it follows that 
$$0=\mu_\sigma^t(x_t)=\mu^{\rho^t}_\sigma(x_t)\ .$$
For $\lambda\in \frak z_L$ we have $\lambda=\lambda_S+\lambda_L$ for 
some $\lambda_S\in\frak s$ and $\lambda_L\in\frak l$ and 
$[\lambda,\lambda_L]=0$. Thus, using the fact that $x_t$ is an $L^\C$-fixed
point, we have 
$$\exp\,is\lambda\cdot x_t=
                            \exp\,is\lambda\cdot \exp\,is\lambda_L\cdot x_t
                           =\exp\,is\lambda_S\cdot x_t\,.$$ 
This implies 
$$\eqalign{
0 &=\mu_{\lambda_S}^t(x_t)   \cr
  &=\mu_{\lambda_S}^{\rho^t}(x_t) \cr
  &=\left({d\over ds}\right)_{s=0}\,
\rho^t(\exp\,is\lambda_S\cdot x_t)\cr
  &=\left({d\over ds}\right)_{s=0}\,
\rho^t(\exp\,is\lambda\cdot x_t) \cr
  &=\mu_\lambda^{\rho^t}(x_t)\ .\cr}$$
Since 
$$\mu^{\rho^t}_\xi(x_t)=
 \mu^{\rho^t}_\tau(x_t)+\mu^{\rho^t}_\lambda(x_t)+\mu^{\rho^t}_\sigma(x_t)
=\mu^{\rho^t}_\tau(x_t)\ ,$$
this implies the claim.
  
\medskip
Now, as we already observed, on $T^\C$ there exists a pluriharmonic 
function $\tilde h^t:T^\C\to\R$ such that 
$\mu^{\tilde h^t }=c^t:=\mu^{\rho^t}(x_t)\in\frak t^*$. Since 
$T^\C=(K^\C/S^\C)/(S^\C L^\C/S^\C)=K^\C/S^\C L^\C$, the natural map 
$q:K^\C /L^\C \to T^\C$ is $K^\C$-equivariant. Thus $h^t:=\tilde h^t\circ q$ 
is a $K$-invariant pluriharmonic function on $X=K^\C/L^\C$ such that 
$\mu^h(x_t)=c^t$. Therefore $\rho^t - h^t$ is a smooth $K$-invariant function 
such that $\omega^t= 2i\partial\bar\partial\,\,(\rho^t- h^t)$ and, since 
$\mu^t(x_t)=\mu^{\rho^t- h^t}(x_t)=0$ and $X$ is connected, 
$\mu^t=\mu^{\rho^t- h^t}$.
\qed

\nbbbigskip
{\bf 7. Action of a torus}

\nbbigskip
Let $T\cong (S^1)^m$ be a torus and $X$ a  complex projective  
manifold with an algebraic action of the complexified torus 
$T^\C\cong (\C^*)^m$. Let $\omega^j$ be $T$-invariant K\"ahler forms
on $X$ with moment maps $\mu^j:X\to\frak t$.  

\smallskip
We say that $\omega^0$ and $\omega^1$ are cohomologous on the
closure $Y$ of a $T^\C$-orbit in $X$ if there is a $T^\C$-equivariant
projective desingularization $p:\tilde Y\to Y$ such that the 
pull back of the forms to $\tilde Y$ are cohomologous, i.e., such that
$p^*\omega^1-p^*\omega^0=2i\partial\bar\partial\, f$ for some $T$-invariant
smooth function $f:X\to \R$.

\smallskip
For $t\in [0,1]$ we set $\omega^t:= (1-t)\,\omega^0 +t\,\omega^1$
and $\mu^t:=(1-t)\,\mu^0 + t\, \mu^1$. Note that $\mu^t$ is a moment 
map with respect to $\omega^t$ and that $\omega^t$ and $\omega^0$ are
cohomologous on the closure of every $T^\C$-orbit in $X$ is this is
the case for $\omega^0$ and $\omega^1$.

\nbbigskip
{\bf Proposition.} {\it If $\omega^0$ and $\omega^1$ are cohomologous
on the closure of every $T^\C$-orbit in $X$, then there is a 
constant $c^t\in\frak t^*$ depending continuously on $t$ such that 
$$X(\mu^0)=X(\mu^t+c^t)\,.$$}

\nbbigskip 
For the proof of the Proposition we consider first the case where 
$T\cong S^1$, i.e., we fix a one dimensional subtorus 
$S^1=\{\exp\,z\xi;\, z\in \R\}$ where $\xi$ is chosen to be a generator 
of the kernel of the one-parameter group $z\to \exp \,z\xi$. With respect 
to this $S^1$-action let $X^{S^1}=\cup F_\alpha$ be the decomposition of the 
set of $S^1$-fixed points  of $X$ into connected  components. 
The set of these components is endowed with a partial order relation 
which is generated by $F_\alpha<F_\beta$. Here we set $F_\alpha<F_\beta$ 
if and only there is a point $x\in X$ such that 
$\buildunder{z\to 0}\lim\, z\cdot x\in F_\alpha$ and 
$\buildunder{z\to\infty }\lim \,z\cdot x\in F_\beta$
where $z\in \C^*=(S^1)^\C$.

\smallskip
Let $\mu^t_\xi:X\to \R$ where $\mu^t_\xi=\langle \mu^t\,,\xi\rangle$ denote 
the moment map with respect to the given $S^1$-action. Since 
$d\mu^t_\xi=\imath_{\xi_X}\omega^t$, the moment map
$\mu^t_\xi$ is constant on every $F_\alpha$.  

\nbbigskip
{\bf Lemma.} {\it If $F_\alpha<F_\beta$, then  
$\mu^0_\xi (F_\alpha)-\mu^0_\xi (F_\beta)=\mu^t_\xi(F_\alpha)-
\mu^t_\xi(F_\beta)$.}

\nbigskip
{\it Proof.} Let $x_0\in X$ be such that 
$\buildunder{z\to 0}\lim\,z\cdot x_0\in F_\alpha$ and 
$\buildunder{z\to\infty }\lim\,z\cdot x_0\in F_\beta$. We may assume that 
the map $\C^*\to \C^*\cdot x_0,\, z\to z\cdot x_0$ is an isomorphism
and extends to a holomorphic map $b: \P_1(\C)\to X$ with $b(0)=x_\alpha$ 
and $b(\infty)=x_\beta$.

\smallskip
Now since by assumption the pull back of 
$\eta:=\omega^t-\omega^0$ to the desingularization
$\P_1(\C)$ of $\overline{\C^*\cdot x_0}$ is cohomologous to zero we have
$$\eqalign{
0 &=\int_{\overline{\C^*\cdot x_0}}\eta     \cr
  &=\int_{\C^*\cdot x_0}\eta                \cr
  &=\int_{\R^+\cdot x_0}\imath_{\xi_X}\eta  \cr
  &=\int_{\R^+\cdot x_0}d(\mu^t_\xi-\mu^0_\xi) \cr
  &=\mu^t_\xi(x_\beta)-\mu^0_\xi(x_\beta)
-(\mu^t_\xi(x_\alpha)-\mu^0_\xi(x_\alpha))\,. \cr}$$
Here $\R^+\cdot x_0$ denotes the $\R^+:=\{z\in \R;\, z>0\}$-orbit 
through $x_0$.

\qed

\nbbigskip
{\it Remark.} Implicitly we used that under the above assumption 
$\omega^0$ and $\omega^1$ are cohomologous on the normalisation of 
$\overline{\C^*\cdot x_0}$.

\nbbigskip
{\it Proof of the Proposition.} The above Lemma implies that there is 
a constant $c^t\in \frak t^*$ depending continuously on $t$
such that $\mu^0$ and $\tilde\mu^t:=\mu^t+c^t$
assume the same values on every component of the set $X^T$ of $T$-fixed
points in $X$.  
Since $\tilde\mu^t(\overline{T^\C\cdot x)})$ is the convex 
hull of the images of $\tilde\mu^t(F_\alpha)$ where 
$F_\alpha\cap \overline{T^\C\cdot x}\ne\emptyset$ (see [A])
it follows that $X(\mu^0)=\{x\in X;\ 0\in \mu^0(\overline{T^\C\cdot x)})\}
=\{x\in X;\ 0\in \tilde\mu^t(\overline{T^\C\cdot x)})\}=X(\tilde\mu^t)$.
\qed

\nbbbigskip
{\bf 8. Action of a semisimple group}

\nbbigskip
Let $G$ be a connected complex semisimple Lie  group with maximal compact
subgroup $K$ and $X$ a projective manifold with an algebraic 
$G$-action.
As in the last section we say that two given closed forms 
$\omega^0$ and $\omega^1$ are cohomologous on the
closure $Y$ of a $G$-orbit in $X$ if there is a $G$-equivariant
desingularization $p:\tilde Y\to Y$ such that 
$p^*\omega^0-p^*\omega^1=2i\partial\bar\partial f$ 
for some smooth function $f:\tilde Y\to \R$.

\nbbigskip
{\bf Proposition.} {\it Let $\omega^j:X\to \frak k^*$, $j=0,1,$ be two 
$K$-invariant K\"ahler forms on $X$ which are cohomologous on every
$G$-orbit closure and let $\mu^j$ be the unique moment map with respect to
$\omega^j$. Then
$$X(\mu^0)=X(\mu^1)\,.$$           }

\nbigskip
{\it  Proof.} 
For $x\in (\mu^0)^{-1}(0)$ set $Y:=\overline{G\cdot x}$ and let 
$p:\tilde Y\to Y$ an equivariant resolution of singularities 
such that $p^*\omega^1 -p^*\omega^0 =2i\partial\bar\partial f$ for a
smooth $K$-invariant function $f$. In particular $f$ is bounded on $G\cdot x$.

\smallskip
Since $G\cdot x$ is closed in $X(\mu^0)$ it follows from the exhaustion
lemma that $\mu^0\vert G\cdot x =\mu^\rho$ for some 
K-invariant plurisubharmonic exhaustion function $\rho:G\cdot x\to \R$.

\smallskip
Therefore $\rho +f$ is likewise an exhaustion and in particular  
has a minimum on $G\cdot x$. Since $\mu^1$ is unique, 
we have $\mu^1=\mu^0+\mu^f=
\mu^{\rho+f}$. Thus $X(\mu^0)\subset X(\mu^1)$ and the reverse 
inclusion follows by symmetry.
\qed

\nbbbigskip
{\bf  9. Reduction to Levi factors}

\nbbigskip
Let $G$ be a connected complex reductive group with 
maximal compact subgroup $K$ and let $X$ be a compact connected 
manifold endowed with a holomorphic action of $G$. 
We assume that there are given $K$-invariant K\"ahler forms
$\omega^j$, $j=0,1,$ on $X$ which are cohomologous on any $G$-orbit
and set 
$$\omega^t=(1-t)\,\omega^0+t\,\omega^1\,, \ t\in [0,1]\,.$$
Moreover, assume that there are $K$-equivariant moment maps
$$\mu^t:X\to \frak k^*\,,\ t\in[0,1]$$
with respect to $\omega^t$ which depend continuously on $t$.
We set $M_K^t:=(\mu^t)^{-1}(0)$. 

\smallskip
Let $Z$ be the center of $K$ and $S$ the semisimple 
part of $K$, i.e.,  $K=Z\cdot S$ where $Z\cap S$ is a finite group
and assume the following condition:
$$X(\mu^t_Z)\ \hbox{is independent of } t\in [0,1]\,. \leqno (*)$$

\nbbigskip
{\bf Lemma 1.} {\it Assume the condition $(*)$ and for $x_0\in X$ let 
$\Omega:=G\cdot x_0$.
Then, for $t\in [0,1]$, 
$$M_K^t\cap \Omega\ne\emptyset$$ 
is an open condition.}

\nbigskip
{\it Proof.} Let $t_0\in [0,1]$ be such that 
$M_K^{t_0}\cap \Omega \ne \emptyset $.
It follows that there exists a smooth curve $\rho^t$ of $K$-invariant
smooth functions so that
$\omega^t=2i\partial {\bar \partial }\,\tau^t $ and
$\mu^t=\mu^{\rho^t}$ on $\Omega$.  
Furthermore, 
since $M_K^{t_0}\cap \Omega \ne \emptyset $,
it follows that $\rho^{t_0}$ is an exhaustion of $\Omega$.  For
$t$ near $t_0$ the function $\rho^t$ has the same convexity properties
as $\rho^{t_0}$ and is therefore likewise an exhaustion (see [H-H], proof
of Lemma 2 in sec. 2).
The points where it has its minimum are those in $M_K^t\cap \Omega$.
\qed

\nbbigskip
{\bf Lemma 2.} {\it Assume the condition $(*)$ and for $x_0\in X$ let 
$\Omega:=G\cdot x_0$.
Then, for $t\in [0,1]$, 
$$M_K^t\cap \Omega\ne\emptyset$$ 
is a closed condition.}

\nbigskip
{\it Proof.} We have to show that 
$M_K^t\cap \Omega \ne \emptyset$ for $t\le t_0$ implies 
$M_K^{t_0}\cap \Omega \ne \emptyset $.

\smallskip
Since $\mu^t_K$ depends continuously on $t$ and $X$ is compact, it follows
that $M_K^{t_0}\cap \overline{G\cdot x_0}\ne \emptyset$.  
Let $y_0\in M_K^{t_0}$.  If   
$G\cdot y_0\ne \Omega $, then by Lemma 1 for $t$ near $t_0$ 
we have that $M_K^t\cap \overline{G\cdot x_0}\subset G\cdot y_0$. 
However the intersection
of $M^t_k$ with $\overline{G\cdot x_0}$ consist of precisely one $K$-orbit,
which would be contrary to $M^t_K\cap \Omega$ also being non-empty.
\qed

\nbbigskip
{\bf Proposition.} {\it Assume that condition $(*)$ is fulfilled.
Then 
$X(\mu^t_K)$ does not depend on $t\in [0,1]$.}

\nbigskip
{\it Proof.} Let $x\in M_K^{t_0}$ and $\Omega :=G\cdot x$. From the above
two Lemma it follows that $M_K^t\cap \Omega \ne \emptyset $
for all $t$. Thus, the condition
that $\Omega $ is a closed $G$-orbit in $X(\mu^t_K)$
is satisfied for some $t$ if and only if this is the case for all $t$.
\qed

\nbbbigskip
{\bf 10. Proof of the Semistability Theorem}

\nbbigskip
For the proof of the Semistability Theorem we need to associate to
a given K\"ahler form one whose cohomology class lies in the 
real span of the ample cone. 
(see [M], \S3).
Let $X$ be a smooth projective variety and denote by $H\in H^2(X,\Q)$
the cohomology class of a hyperplane section.

\smallskip
Let ${\cal C}_1$ 
be the subspace of the second rational homology group
$\H_2(X,\Q)$ which is spanned by the images of closed analytic curves and
${\cal C}_{n-1}$ 
the subspace of $H^2(X,\Q)$ spanned by divisors, or, 
what is the same, Chern classes of holomorphic line bundles.

\nbbigskip
{\bf Lemma 1.} {\it The pairing ${\cal C}_1 \times {\cal C}_{n-1} 
\longrightarrow \R$ which is induced by associating to a line bundle 
$L$ and a curve $C$ the intersection number $L\cdot C:=\hbox{\rm deg}\, L_C$ 
is perfect.}

\nbigskip
{\it Proof.} The result is well known for surfaces.
Let $S$ be a generic intersection of $n-2$ hyperplane sections.
If a divisor $D$ is numerically trivial its restriction $D_S$ to $S$
is numerically trivial as well and therefore its cohomology class in 
$H^2(S,\Q)$ is zero. By the Lefschetz hyperplane sections theorem 
the restriction map $H^2(X,\Q)\longrightarrow H^2(S,\Q)$ is injective 
and therefore the cohomology class of $D$ is zero.
Let now $C$ be a numerically trivial 1-cycle on $X$. By the Hard 
Lefschetz theorem
there is a divisor $D$ such that $C$ is homologous to
$H^{n-2}D$. If $\gamma=H^{n-2}\Gamma$ is another 1-cycle 
$D \cdot \gamma= D \cdot H^{n-2} \cdot \Gamma=C \cdot \Gamma=0$ so that
$D$ is numerically trivial and homologous to zero, and
therefore  $C$ is homologous to zero.
\qed

\nbbigskip
{\bf Lemma 2.} {\it Let $\omega$ be a K\"ahler form on $X$. 
Then there exist a K\"ahler form $\widetilde \omega$ whose cohomology class
$[\widetilde{\omega }]$ lies in the span of the ample cone  such 
that 
$$\int_C\,\widetilde\omega =\int_C\, \omega$$  
holds for all one-dimensional analytic cycles $C$.
}

\nbigskip
{\it Proof.} 
Consider the linear map 
$\lambda:{\cal C}_1\to \R, \, \lambda(C)=\int_C\omega,$ 
given by $\lambda(C)=\int_C\omega$. By Lemma 1 there is a class
$\tilde D$ in ${\cal C}_{n-1}$ such that 
$\lambda(C)=\tilde D \cdot C$ for all 1-cycles $C$.
Since $\tilde D$ is a divisor, the cohomology class of $\tilde D$ lies
in the span of the ample cone.

\smallskip
As a consequence of the
Kleiman ampleness criterion 
the cohomology class of $\tilde D$ contains
a K\"ahler form
$\widetilde \omega$.

\qed

\nbbigskip
{\bf Lemma 3.} {\it Let $Y$ be a connected smooth projective variety and 
assume that $G$ has an open orbit on $Y$. Then there are no non-zero 
holomorphic $p$-forms on $Y$ for $p\ge 1$.}

\nbigskip
{\it Proof.} For $\xi\in \frak k$ let $\xi_Y$ denote the corresponding 
holomorphic vector field on $Y$. The complex one-parameter subgroup of
$G$ which is generated by $\xi$ has a fixed point on $Y$. 
This is seen most easily by considering the closure of this group in $G$ 
which is an algebraic torus and then an induction argument shows that 
every algebraic torus action on a projective variety has a 
fixed point.

\smallskip
Now let $\alpha$ be a holomorphic $p$-form on $Y$ and $z_0$ a point in the
open $G$-orbit. For $\xi_j\in \frak k,\, j=1,\ldots,p$ the 
holomorphic function $\alpha(\xi_{1Y},\ldots,\xi_{pY})$ is constant on $Y$ 
and therefore identically zero.
Since $\{\xi_Y(z_0);\ \xi\in \frak k\}$ span the complex tangent space
at $z_0$, it follows that $\alpha$ is zero on the open $G$-orbit. Hence
$\alpha$ is identically zero on $Y$. 
\qed

\nbbigskip

Under the same assumption on $Y$ as in the above Lemma we have the following

\nbbigskip
{\bf Corollary.} {\it If $\alpha$ is a smooth closed $(1,1)$-form on $Y$
such that $\int_C\alpha=0$ for every one-dimensional analytic cycle $C$,
then $\alpha=2i\partial\bar\partial\, f$ for some smooth function $f:Y\to\R$.}

\nbigskip
{\it Proof.} This follows from Lemma 3 since 
$\H^{1,1}(Y)=\H^2(Y,\C)=\H^2(Y,\Z)\otimes\C$ holds (see the proof of Lemma 2).
\qed

\nbbigskip
{\bf Lemma 4.} {\it If the $K$-action on $X$ is Hamiltonian with
respect to the $K$-invariant K\"ahler form $\omega$, then it is also
Hamiltonian with respect to any other $K$-invariant K\"ahler 
form $\tilde\omega$.}

\nbigskip
{\it Proof.} Since $X$ is algebraic and the $K^\C$-action on $X$ is
assumed to be algebraic, every one-parameter subgroup $t\to \exp\, t\xi$,
$t\in \R$, $\xi\in \frak k$, has a fixed point in $X$. Using the 
Hodge decomposition (see [F]) one sees that to every
$\xi\in \frak k$ there exist a corresponding Hamiltonian function 
$\mu_\xi:X\to \R$ such that 
$$d\mu_\xi=\imath_{\xi_X}\omega\,.$$
Moreover, $X$ is compact and therefore we may assume that
$$\int_X\mu_\xi\,\omega^n=0$$
for every $\xi$. Now it follows that 
$\int_X\{\mu_{\xi_1},\mu_{\xi_2}\}\,\omega^n=0$ for all $\xi_1,\xi_2\in \frak k$
and this implies that $\{\mu_{\xi_1},\mu_{\xi_2}\}=\mu_{[\xi_1,\xi_2]}$, i.e.,
since $K$ is connected $\mu:X\to\frak k,\, \mu(x)(\xi):=\mu_\xi(x)$
is $K$-equivariant. Here $\{\ ,\ \}$ denotes the Poisson brackets on
${\cal C}^\infty(X)$ given by $\omega$.
\qed

\nbbigskip
{\it Proof of the Semistability Theorem.}
Given a smooth  $K$-invariant K\"ahler form $\omega$ on a smooth 
projective $G$-variety $X$ we already know from Lemma 2 that there is a K\"ahler form 
$\tilde \omega$ on $X$ which lies in the $\R$-span of the ample cone
of $X$ such that 
$$\int_C\omega=\int_C\tilde\omega \leqno (*)$$ 
on every analytic curve $C$ in $X$. Since $K$ is assumed to be connected 
the cohomology class of $\tilde\omega$ is $K$-invariant. 
Hence, after integration 
over the compact group $K$, we may assume that 
$\tilde\omega$ is $K$ invariant and still satisfies $(*)$.   

\smallskip
Now it follows from Lemma 4, the Proposition in sec. 7 and the existence of
a moment map in the semisimple case that there is a moment map
$\mu^t:X\to\frak k^*$ with 
respect to $\omega^t$ where $\omega^t=(1-t)\,\tilde\omega+ t\, \omega$
such that the $\mu^t$ depends continuously on $t$ and such that 
$X(\mu^t_Z)=X(\mu^0_Z)$ for all $t\in [0,1]$. 
Here $Z$ denotes the connected component of the center of $K$.

\smallskip
Moreover $(*)$ implies 
$$\int_C\omega^t=\int_C\omega \,.\leqno (**)$$ 
Since the closure of every $G$-orbit in $X$ has an equivariant 
algebraic desingularization, it follows from the above Corollary
than the forms are cohomologous on the closure of every $G$-orbit.
The statement of the theorem now follows from the Proposition in sec. 9
and the Proposition in sec. 5.
\qed

\nbbbigskip
\centerline{\it References}

\nbbigskip

{\parindent=1,8cm
\litem{[A]} Atiyah, M. F.: {\it Convexity and commuting Hamiltonians.}
           Bull. London Math. Soc. 14, 1--15 (1982)
\smallskip

\litem{[A-L]} Azad, H.; Loeb, J. J.: {\it Pluri-subharmonic functions and the 
             Kempf-Ness theorem.}
             Bull. London Math. Soc. 25, 162--168 (1993)
\smallskip

\litem{[Ch]} Chevalley, C.: Theory of Lie groups, I.
             Princeton, Princeton University Press, 1946.

\smallskip
\litem{[D-W]} Dolgachev, I.; Hu, Y.: {\it Variation of geometric in\-va\-riant
              theory quotients.}  Pre\-print (1994)

\smallskip
\litem{[F]} Frankel, T. T.: {\it Fixed points and torsions on K\"ahler 
            manifolds.}
            Ann. of Math. 70, 1--8 (1959) 
\smallskip
\litem{[G-S]} Guillemin,V.; Sternberg, S.: Symplectic techniques in physics.
              Cambridge University Press 1984 

\smallskip
\litem{[H]} Heinzner, P.: {\it Geometric invariant theory on Stein spaces.}
            Math. Ann. 289, 631--662 (1991)

\smallskip
\litem{[H-H]} Heinzner, P.; Huckleberry, A. T.: {\it K\"ahlerian 
              potentials and convexity properties of the moment map.} 
              Invent. math. 126 , 65--84 (1996)

\smallskip
\litem{[H-H-L]} Heinzner, P.; Huckleberry, A. T.; Loose, F.: {\it K\"ahlerian
extensions of the symplectic reduction.}  
J. reine und angew. Math. 455, 123--140  (1994)

\smallskip
\litem{[H-M-P]} Heinzner, P.; Migliorini, L.; Polito, M.: 
{\it Semistable quotients.}  
Annali della Scuola Normale Superiore di Pisa (1997) (to appear)
 
\smallskip
\litem{[H-L]} Heinzner, P.; Loose, F.: {\it Reductions of complex Hamiltonian 
              $G$-spaces.}
              Geometric and Functional Analysis 4, 288--297 (1994)

\smallskip
\litem{[K-N]} Kempf, G.; Ness, L.: {\it The length of vectors in representation
             spaces.}
             In: Lect. Notes Math. vol. 732, Springer--Verlag, Berlin 
             Heidelberg New York, 233--243, 1979

\smallskip
\litem{[K]} Kirwan, F.: Cohomology of quotients in symplectic and
            algebraic geometry.
            Mathematical notes 31, Princeton University Press, Princeton
            New Jersey, 1984

\smallskip
\litem{[M]} Moishezon, B. G.: {\it On n-dimensional compact complex varieties
            with n algebraically independent meromorphic functions.}
            Amer. Math. Soc. Transl. (2) 63, 51--177, (1967)
\smallskip
\litem{[M-F-K]} Mumford, D.; Fogarty, J.; Kirwan, F.:
                Geometric Invariant Theory.
                Lecture Notes in Math. Springer 3rd ed. Ergeb.Math.
                Berlin Heidelberg New York 1994

\smallskip
\litem{[N]} Ness, L.: {\it A stratification of the null cone via the moment 
            map.}
            Amer. J. Math. 106, 1326--1329 (1983)
\smallskip

\litem{[S]} Sjamaar, R.: {\it Holomorphic slices, symplectic reduction and
            multiplicities of representations.}
            Ann. of Math. 141, 87--129 (1995)

\smallskip
\litem{[T]}  Thaddeus, M.:  {\it Geometric invariant theory and flips.} 
             Preprint (1994)

}

 \bigskip 
 Peter Heinzner

 Brandeis University

 Waltham MA 02254-9110

 USA

 e-mail: heinzner@max.math.brandeis.edu

 \bigskip
 Luca Migliorini

 Dipartimento di Matematica Applicata "G.Sansone"

 via  S. Marta 3 

 I-50139 Firenze 

 Italy
  
 e-mail: lucap@ingfi1.ing.unifi.it

\end